\input amstex
\documentstyle{amsppt}
\magnification=\magstephalf
\pagewidth{6.25truein}
\pageheight{9.25truein}
\parskip=10pt
\voffset=0pt
\font\smallit=cmti10
\font\smalltt=cmtt10
\font\smallrm=cmr8
\NoBlackBoxes
\nologo
\pageno=1
\TagsOnRight
\def\Z{\Bbb Z}
\def\N{\Bbb N}

\def\bg{\bigg}
\def\({\bg(}
\def\[{\bg[}
\def\){\bg)}
\def\]{\bg]}
\def\t{\text}
\def\f{\frac}
\def\mo{\roman{mod}}

\def\eq{\equiv}

\def\ls{\leqslant}
\def\gs{\geqslant}

\def\ve{\varepsilon}

\topmatter \leftheadtext{\hskip 8pt \smalltt INTEGERS: \smallrm
Electronic Journal of Combinatorial Number Theory \smalltt 7
(2007), \#A56} \rightheadtext{\smalltt INTEGERS: \smallrm
Electronic Journal of Combinatorial Number Theory \smalltt 7
(2007), \#A56 \hskip 8pt}
\endtopmatter

\document
\vskip -40pt \centerline{\hskip 8pt \smalltt INTEGERS: \smallrm
ELECTRONIC JOURNAL OF COMBINATORIAL NUMBER THEORY \smalltt 7
(2007), \#A56 }
\vskip 40pt
\centerline{\bf MIXED SUMS OF SQUARES
AND TRIANGULAR NUMBERS (II)}
\vskip 20pt

\centerline{\smallit Song Guo, Department of  Mathematics, Huaiyin
Teachers College, Huaian 223001, P. R. China} \centerline{\tt
guosong77\@sohu.com}\vskip 20pt \centerline{\smallit Hao Pan,
Department of Mathematics, Shanghai Jiaotong University, Shanghai
200240, P. R. China} \centerline{\tt haopan79\@yahoo.com.cn}
\vskip 20pt \centerline{\smallit Zhi-Wei Sun\footnote{This author
is responsible for communications, and supported by the National
Science Fund (Grant No. 10425103) for Distinguished Young Scholars
in China.  Homepage: {\tt http:/\!/math.nju.edu.cn/$\sim$zwsun}},
Department of Mathematics, Nanjing University, Nanjing 210093, P.
R. China} \centerline{\tt zwsun\@nju.edu.cn} \vskip 30pt

\centerline{\smallit Received: 4/16/07, Revised: 11/12/07,
Accepted: 12/6/07, Published: 12/20/07} \vskip 30pt

\centerline{\bf Abstract} For $x\in\Z$ let $t_x$ denote the
triangular number $x(x+1)/2$. Following a recent work of Z. W.
Sun, we show that every natural number can be written in any of
the following forms with $x,y,z\in\Z$:
$$x^2+3y^2+t_z,\ x^2+3t_y+t_z,\ x^2+6t_y+t_z,\ 3x^2+2t_y+t_z,\ 4x^2+2t_y+t_z.$$
This confirms a conjecture of Sun.

\noindent

\baselineskip=15pt \vskip 30pt

\subhead \nofrills{\bf 1. Introduction}\endsubhead

In 1916 S. Ramanujan [6] found all those positive integers $a,b,c,d$
such that every natural number can be written in the form $ax^2+by^2+cz^2+dw^2$
with $x,y,z,w\in\Z$.

Let $a,b,c$ be positive integers with $a\ls b\ls c$.
In 2005 L. Panaitopol [5] showed that
any positive odd integer can be written as $ax^2+by^2+cz^2$ with $x,y,z\in\Z$, if and only if
the vector $(a,b,c)$ is $(1,1,2)$ or $(1,2,3)$ or $(1,2,4)$.

As usual, for any $x\in\Z$ we call $t_x=x(x+1)/2$ a {\it
triangular number}. In 1862 J. Liouville (cf. L. E. Dickson [1,
p.\,23]) determined those positive integers $a,b,c$ for which any
natural number can be written as $at_x+bt_y+ct_z$ with
$x,y,z\in\Z$.

Let $n\in\N=\{0,1,2,\ldots\}$. As observed by L. Euler (cf. [1,
p.\,11]), the fact that $8n+1$ is a sum of three squares (of integers) implies
that $n$ can be expressed as a sum of two squares and a triangular
number. According to [1, p.\,24], E. Lionnet stated, and V. A.
Lebesgue [3] and M. S. R\'ealis [7] showed that $n$ is also a sum
of two triangular numbers and a square. In 2006 this was re-proved
by H. M. Farkas [2] via the theory of theta functions.

In [8] Z. W. Sun investigated mixed sums of squares and triangular
numbers systematically, and he mainly proved the following result.

\proclaim{Theorem 1 {\rm (Sun [8])}} {\rm (i)} Any natural number
is a sum of an even square and two triangular numbers,
and each positive integer is a sum of a triangular number
plus $x^2+y^2$ for some $x,y\in\Z$ with
$x\not\eq y\ (\mo\ 2)$ or $x=y>0$.

{\rm (ii)} Let $a,b,c$ be positive integers with $a\ls b$. If
every $n\in\N$ can be written as $ax^2+by^2+ct_z$ with
$x,y,z\in\Z$, then $(a,b,c)$ is among the following vectors:
$$\gather(1,1,1),\ (1,1,2),\ (1,2,1),\ (1,2,2),\ (1,2,4),
\\(1,3,1),\ (1,4,1),\ (1,4,2),\ (1,8,1),\ (2,2,1).
\endgather$$

{\rm (iii)} Let $a,b,c$ be positive integers with $b\gs
c$. If every $n\in\N$ can be written as $ax^2+bt_y+ct_z$
with $x,y,z\in\Z$, then $(a,b,c)$ is among the following vectors:
$$\gather(1,1,1),\ (1,2,1),\ (1,2,2),\ (1,3,1),\ (1,4,1),\ (1,4,2),\ (1,5,2),
\\ (1,6,1),\ (1,8,1),
\ (2,1,1),\ (2,2,1),\ (2,4,1),\ (3,2,1),
\ (4,1,1),\ (4,2,1).
\endgather$$
\endproclaim

Sun also reduced the converses of (ii) and (iii) to
Conjectures 1 and 2 of [8]. In this paper we
prove his second conjecture, namely we establish the following theorem.

\proclaim{Theorem 2} Every $n\in\N$ can be expressed in any of the following forms with $x,y,z\in\Z$:
$$x^2+3y^2+t_z,\ x^2+3t_y+t_z,\ x^2+6t_y+t_z,\ 3x^2+2t_y+t_z,\ 4x^2+2t_y+t_z.$$
\endproclaim

\vskip 30pt

\subhead \nofrills{\bf 2. Proof of Theorem 2}\endsubhead

The following theorem is well-known (cf. [4, pp.\,17-23]).

 \proclaim{Gauss-Legendre Theorem} A natural number can be written as a sum of three squares
 of integers
 if and only if it is not of the form $4^k(8l+7)$ with $k,l\in\N$.
 \endproclaim

We also need an identity of Jacobi which can be verified directly.

\proclaim{Jacobi's Identity} We have
 $$3(x^2+y^2+z^2)=(x+y+z)^2+2\(\f{x+y-2z}{2}\)^2+6\(\f{x-y}{2}\)^2.$$
 \endproclaim

 \noindent{\it Proof of Theorem 2}.
(i) By the Gauss-Legendre theorem, $8n+3=x^2+y^2+z^2$ for some
$x,y,z\in\Z$. Clearly, each of $x,y,z$ is congruent to $1$ or $-1$
modulo $4$. Without any loss of generality, we simply let $x\eq
y\eq z\eq1\ (\mo\ 4)$. Two of $x,y,z$ are congruent modulo $8$,
say, $x\eq y\ (\mo\ 8)$. Set
$$x_0=\f{x-y}8,\ \ y_0=\f{x+y-2}4\ \ \t{and}\ \ z_0=\f{z-1}2.$$
Then
$$8n+3=2\(\f{x-y}2\)^2+2\(\f{x+y}2\)^2+z^2=2(4x_0)^2+2(2y_0+1)^2+(2z_0+1)^2$$
and hence $n=4x_0^2+2t_{y_0}+t_{z_0}$.
\medskip

(ii) By the Gauss-Legendre theorem,
$12(4n+2)=x^2+y^2+z^2$ for some $x,y,z\in \Z$. As
$x^2+y^2+z^2\eq0\ (\mo\ 3)$, we can choose suitable
$\ve_1,\ve_2,\ve_3 \in \{\pm1\}$ such that
 $\ve_1 x\eq \ve_2y \eq \ve_3z \eq 0 \ \t{or} \ 1 \ (\mo \ 3)$.
Therefore we may simply let $x\eq y \eq z \ (\mo \ 3)$.
Since $x^2+y^2+z^2\eq 8 \ (\mo \ 16)$,  $x,y,z$ are all even and exactly one of them
is divisible by 4. Suppose that $x \eq y+2 \eq z+2 \eq 0\ (\mo \ 4)$.
It is easy to see that
$$x+y+z\eq 0\ (\mo \ 12),\ x+y-2z\eq 6\ (\mo \ 12),\ x-y\eq 6\ (\mo \ 12).$$

Set
$$x_0=\f{x+y+z}{12},\ \ y_0=\f{x+y-2z-6}{12}\ \ \t{and}\ \ z_0=\f{x-y-6}{12}.$$
By Jacobi's identity,
$$\aligned 36(4n+2)=&3(x^2+y^2+z^2)
\\=&(12x_0)^2+2(6y_0+3)^2+6(6z_0+3)^2 \\
=&144x_0^2+72y_0(y_0+1)+18+216z_0(z_0+1)+54.
\endaligned$$
It follows that $n=x_0^2+t_{y_0}+3t_{z_0}$.
\medskip

(iii) Let $\ve\in\{0,1,3\}$. By the Gauss-Legendre theorem,
$24n+3+6\ve=x^2+y^2+z^2$ for some $x,y,z\in \Z$. As $3\mid
x^2+y^2+z^2$, without loss of generality we may assume that $x\eq
y\eq z\eq0\ \t{or}\ 1\ (\mo\ 3).$ Applying Jacobi's identity, we
obtain that
$$72n+9+18\ve =(x+y+z)^2+2\(\f{x+y-2z}{2}\)^2+6\(\f{x-y}{2}\)^2.$$
Recall that $x^2+y^2+z^2\eq 3+6\ve \ (\mo \ 8)$.
If $\ve=0$, then $x,y,z$ are odd, and two of them are congruent modulo $4$,
say, $x\eq y\ (\mo\ 4)$.
In the case $\ve=1$, we may suppose that $x\eq y\eq z-1\eq0\ (\mo\ 2)$
and $x\eq y\ (\mo\ 4)$. When $\ve=3$, we may assume that $x\eq2\ (\mo\ 4)$, $4\mid y$ and $2\nmid z$.
Clearly,
$$x+y+z\eq 3\ (\mo\ 6),\ \ \ x+y-2z\eq \cases0 \ (\mo \ 12)&\t{if} \ \ve=0,3,\\
6 \ (\mo \ 12)&\t{if} \ \ve=1,\endcases$$
and $$x-y\eq \cases
0 \ (\mo \ 12) &\t{if}\ \ve=0,1,\\
6 \ (\mo \ 12) &\t{if}\ \ve=3.
\endcases$$

Set
$$ \align x_0=& \cases
(x+y-2z)/12 &\t{if}\ \ve=0,\\
(x-y)/12 &\t{if}\ \ve=1,\\
(x+y-2z)/12 &\t{if}\ \ve=3,
\endcases
\\y_0=&\cases
(x-y)/12 &\t{if}\ \ve=0,\\
(x+y-2z-6)/12 &\t{if} \ \ve=1,\\
(x-y-6)/12 &\t{if}\ \ve=3,
\endcases
\endalign$$
and $z_0=(x+y+z-3)/6$. By the above,
$$\align 72n+9+18\ve=&\cases (6z_0+3)^2+2(6x_0)^2+6(6y_0)^2&\t{if}\ \ve=0,
\\(6z_0+3)^2+2(6y_0+3)^2+6(6x_0)^2&\t{if}\ \ve=1,
\\(6z_0+3)^2+2(6x_0)^2+6(6y_0+3)^2&\t{if}\ \ve=3,\endcases
\\=&\cases
72 x_0^2+216y_0^2+36z_0(z_0+1)+9 &\t{if}\ \ve=0,\\
216x_0^2+72y_0(y_0+1)+36z_0(z_0+1)+27 &\t{if}\ \ve=1,\\
72 x_0^2+216y_0(y_0+1)+36z_0(z_0+1)+63 &\t{if}\ \ve=3.
\endcases
\endalign$$
It follows that
$$n=\cases
x_0^2+3y_0^2+t_{z_0}&\t{if}\ \ve=0,\\
3x_0^2+2t_{y_0}+t_{z_0}&\t{if}\ \ve=1,\\
x_0^2+6t_{y_0}+t_{z_0} &\t{if}\ \ve=3.
\endcases$$

\medskip

Combining (i)--(iii) we have completed our proof of Theorem 2. \qed

\vskip 30pt

\subhead \nofrills{\bf References}\endsubhead

\ref\no1\by L. E. Dickson\book
History of the Theory of Numbers, {\rm Vol. II} \publ AMS Chelsea
Publ., 1999\endref

\ref\no 2\by H. M. Farkas\paper Sums of squares and triangular
numbers \jour Online J. Anal. Combin.\vol 1\yr 2006\pages \#1, 11
pp. (electronic)\endref

\ref\no3\by V. A. Lebesgue\paper Questions 1059,1060,1061
(Lionnet) \jour Nouv. Ann. Math.\vol 11\yr 1872\pages
516--519\endref

\ref\no4\by M. B. Nathanson\paper Additive
Number Theory: the Classical Bases\publ Graduate Texts in Math. 164, Springer, New York, 1996\endref

\ref\no5\by L. Panaitopol \paper On the representation of natural numbers as sums
of squares \jour Amer. Math. Monthly\vol 102\yr 2005\pages
168--171\endref

\ref\no6\by S. Ramanujan\paper On the expression of a number in
the form $ax^2+by^2+cz^2+du^2$\jour Proc. Camb. Philo. Soc.\vol 19\yr 1916\pages 11--21\endref

\ref\no7\by M. S. R\'ealis\paper Scolies pour un th\'eoreme
d'arithm\'etique \jour Nouv. Ann. Math.\vol 12\yr 1873\pages
212--217\endref

\ref\no8\by Z. W. Sun\paper Mixed sums of
squares and triangular numbers \jour Acta Arith. \vol 127\yr 2007\pages 103--113\endref

\enddocument